\documentclass[12pts]{svmult}

\usepackage{hyperref}

\usepackage{mathptmx}       
\usepackage{helvet}         
\usepackage{courier}        
\usepackage{type1cm}        

\usepackage{makeidx}         
\usepackage{graphicx}        
\usepackage{multicol}        
\usepackage[bottom]{footmisc}

\usepackage{amsmath}
\usepackage{amssymb}
\usepackage{amscd}
\usepackage{indentfirst}

\newcommand{\EE}{\mathbb{E}}

\newcommand{\C}{\mathbb{C}}

\newtheorem{proph1}{Proposition}
\newtheorem{thmh1}{Theorem}
\newtheorem{proph2}[proph1]{Proposition}
\newtheorem{lemh1}{Lemma}
\newtheorem{lemh2}[lemh1]{Lemma}
\newtheorem{lemh3}[lemh1]{Lemma}
\newtheorem{conjecture1}{Conjecture}
\newtheorem{proph3}[proph1]{Proposition}
\newtheorem{lemh4}[lemh1]{Lemma}

\makeindex             
\begin{document}

\title*{A note on Helson's conjecture on moments of random multiplicative functions}
\author{Adam J. Harper\thanks{The first author is supported by a research fellowship at Jesus College, Cambridge.}, Ashkan Nikeghbali, and Maksym Radziwi\l\l}
%
%
%
%

\institute{Adam J. Harper \at Jesus College, Cambridge CB5 8BL, England, \email{A.J.Harper@dpmms.cam.ac.uk}
\and Ashkan Nikeghbali \at University of Z\"{u}rich, Institute of Mathematics, Winterthurerstrasse 190, CH-8057 Z\"{u}rich, \email{ashkan.nikeghbali@math.uzh.ch}
\and Maksym Radziwi\l\l \at Department of Mathematics, Rutgers University, Hill Center for the Mathematical Sciences,
110 Frelinghuysen Rd., Piscataway, NJ 08854-8019, \email{maksym.radziwill@gmail.com}}

\maketitle


\begin{flushright}
\textit{To Prof. Helmut Maier on the occasion of his sixtieth birthday}
\end{flushright}
\vspace{10mm}



\section{Introduction}

In this note we are interested in cancellations in sums of multiplicative functions. It is well known that
$$
M(x) := \sum_{n \leq x} \mu(n) = O(x^{1/2 + \varepsilon})
$$
is equivalent to the Riemann Hypothesis. On the other hand
it is also a classical result 
that $M(x) > x^{1/2 - \varepsilon}$ for a sequence of arbitrarily large $x$.
It is in fact conjectured that
$$
\overline{\underline{\lim_{x \rightarrow \infty}}} \  \frac{M(x)}{\sqrt{x} (\log\log\log x)^{\tfrac 54}} = \pm B
$$
for some constant $B > 0$ (see \cite{ng}).

Wintner \cite{wintner} initiated the study of what happens for a generic multiplicative
function which is as likely to be $1$ or $-1$ on the primes. Consider $f(p)$, a sequence of independent random variables taking values $\pm 1$ with probability $1/2$ each (i.e. Rademacher random variables), and define a multiplicative function supported on squarefree integers $n$ by
$$ f(n) := \prod_{p | n} f(p) . $$
We shall refer to such a function as a {\em Rademacher random multiplicative function}. By the three series theorem, the Euler product
$F(s) := \prod_{p} ( 1+ f(p) p^{-s} )$ converges
almost surely for $\Re s > \tfrac 12$. From this Wintner deduced 
that
$$
\sum_{n \leq x} f(n) \ll x^{1/2 + \varepsilon}\ \ \text{almost surely (a.s.)}
$$
Since then the problem of the behavior of $\sum_{n \leq x} f(n)$ has attracted considerable attention \cite{chatsound, halasz, harpergp, harperlimits, hough, tenenbaum}. A closely related model is to let $f(p)$ be uniformly distributed on the complex unit circle (i.e. Steinhaus random variables), and then define $f(n) := \prod_{p^{\alpha} || n} f(p)^{\alpha}$ for all $n$. We shall refer to such a function as a {\em Steinhaus random multiplicative function}.

Very recently mean values of random multiplicative functions arose 
in connection with harmonic analysis. In his last paper Helson \cite{helson}
considered the question of generalizing Nehari's theorem to the
infinite polydisk. He noticed that the generalization could be
disproved if one could show that
\begin{equation} \label{helson}
\lim_{T \rightarrow \infty} \frac{1}{T} \int_{0}^{T} \Big | \sum_{n \leq N} n^{-it} \Big | dt = o(\sqrt{N}).
\end{equation}
Using Bohr's identification, we have
\begin{equation} \label{equivalent}
\left( \mathbb{E} \Bigg | \sum_{n \leq N} f(n) \Bigg |^{2q} \right)^{1/2q} = \left( \lim_{T \rightarrow \infty} \frac{1}{T}
\int_{0}^{T} \Bigg | \sum_{n \leq N} n^{-it} \Bigg |^{2q} dt \right)^{1/2q}
\end{equation}
for all $2q > 0$, and with $f(n)$ a Steinhaus random multiplicative function. 
Therefore (\ref{helson}) is equivalent to
\begin{equation} \label{equivalent2}
\mathbb{E} \Big | \sum_{n \leq N} f(n) \Big | = o(\sqrt{N}) ,
\end{equation}
with $f(n)$ a Steinhaus random multiplicative function.
Helson justified his belief in (\ref{helson}) by observing that $N(it) := \sum_{n \leq N} n^{-it}$ is the
multiplicative analogue of the classical Dirichlet kernel $D(\theta) := \sum_{|n| \leq N} e^{2\pi i n \theta}$. 
Since $\| D \|_1 = o(\| D\|_2)$ Helson conjectured that the same phenomenon should happen
for the ``multiplicative analogue'' $N(it)$. Another reason one might believe the large cancellation in (\ref{helson}) to be possible is that on the $\tfrac 12$-line one has
$$ \lim_{T \rightarrow \infty} \frac{1}{T} \int_{0}^{T} \Big | \sum_{n \leq N} \frac{1}{n^{1/2+it}} \Big | dt \ll \log^{1/4+o(1)} N , $$
as follows from the work of Bondarenko, Heap and Seip~\cite{bondhs}. This bound is stronger than one would expect assuming only squareroot cancellation, which would suggest a size more like $\log^{1/2} N$.

Recently Orteg\`{a}-Cerda and Seip \cite{OrtegaSeip} proved that Nehari's theorem doesn't extend to 
the infinite polydisk. However
the problem of establishing (\ref{helson}) remained. There are now reasons to
believe that (\ref{helson}) is false. In a recent paper Bondarenko and Seip \cite{BondarenkoSeip} showed
that the first absolute moment is at least $\sqrt{N} (\log N)^{-\delta + o(1)}$
for some small $\delta < 1$. 
Our primary goal in this note is to improve further on the lower bounds for (\ref{equivalent}).
Our results also work for Rademacher random multiplicative functions.
\begin{theorem} \label{thm:main}
Let $f(n)$ be a Rademacher or Steinhaus random multiplicative function. Then, 
$$
\mathbb{E} \Bigg | \sum_{n \leq N} f(n) \Bigg | \gg \sqrt{N} (\log\log N)^{-3 + o(1)}
$$
as $N \rightarrow \infty$.
\end{theorem}
The main input in the proof of Theorem \ref{thm:main} is the work~\cite{harpergp} of the first named
author on lower bounds for sums of random multiplicative functions. 
Using H\"{o}lder's inequality, we can extend the result of Theorem \ref{thm:main} to $L^q$ norms. 
\begin{theorem} \label{thm:lpmain}
Let $f(n)$ be a Rademacher or Steinhaus random multiplicative function and
let $0 \leq q \leq 1$. Then, 
$$
\mathbb{E} \Bigg | \sum_{n \leq N} f(n) \Bigg |^{2q} \gg N^q (\log\log N)^{-6 + o(1)}.
$$
\end{theorem}
Theorem \ref{thm:main} and Theorem \ref{thm:lpmain} suggest it is rather unlikely that Helson's conjecture is true. See Conjecture \ref{theconjecture}, below.

In addition to the above results, we establish an asymptotic estimate for the $2k$-th moment 
when $k$ 
is a positive integer. 
\begin{theorem} \label{thm:asymp}
Let $k \in \mathbb{N}$.  
Suppose that $f(n)$ is a Steinhaus random multiplicative function.
 Then, as $N \rightarrow \infty$, 
$$
\mathbb{E} \Bigg | \sum_{n \leq N} f(n) \Bigg |^{2k} \sim \binom{2k-2}{k-1} k^{-(k-1)} \cdot c_k \gamma_k \cdot N^k \cdot (\log N)^{(k-1)^2}.
$$
where
$\gamma_k$ is the volume of Birkhoff polytope $\mathcal{B}_{k}$, defined as
the $(k-1)^2$ dimensional volume of the set of $(u_{i,j}) \in \mathbb{R}_{+}^{k^2}$ such that,
\begin{align*}
& \text{for each } i \leq k:
\sum_{1 \leq j \leq k} u_{i,j} = 1 \\\text{ and } & \text{for each } j \leq k:
\sum_{1 \leq i \leq k} u_{i,j} = 1 ,
\end{align*}
and 
$$
c_k = \prod_{p} \Big (1 - \frac{1}{p} \Big)^{k^2} \cdot \Big (1  + \sum_{\alpha \geq 1} \frac{\binom{\alpha + k - 1}{k-1}^2}{p^{\alpha}} \Big ).
$$
\end{theorem}
Note that $\mathcal{B}_k$ is a $(k-1)^2$ dimensional object embedded in a $k^2$ dimensional space. The $(k-1)^2$ dimensional volume of $\mathcal{B}_k$ is equal (see e.g. section 2 of Chan and Robbins~\cite{Robbins}) to $k^{k-1}$ times the full-dimensional volume of the set of $(u_{i,j})_{i,j \leq k-1} \in \mathbb{R}^{(k-1)^2}$ such that, for all $i,j \leq k-1$, 
\begin{align*}
\sum_{j \leq k-1} u_{i,j} \leq 1 \text{ and } \sum_{i \leq k-1} u_{i,j} \leq 1 \text{ and } \sum_{i,j \leq k-1} u_{i,j} \geq k-2.
\end{align*}
The latter is how the volume of $\mathcal{B}_{k}$ will actually arise in our calculations. 

It is worth pointing out that finding a closed formula for the volume of the 
Birkhoff polytope $\mathcal{B}_k$ is a notorious open question and would be 
of interest in enumerative combinatorics, statistics and computational geometry (see \cite{IgorPak}). 
There are evaluations of $\text{Vol}(\mathcal{B}_k)$
for small values of $k$, (see \cite{beckpixton} and \cite{Robbins}), 
$$
\text{Vol}(\mathcal{B}_3) = \frac{3 \cdot 3^2}{2^2!} \ , \ \text{Vol}
(\mathcal{B}_4) = \frac{352 \cdot 4^3}{3^2!} \ , \ \text{Vol}(\mathcal{B}_5) = 
\frac{4718075 \cdot 5^4}{4^2!}, \ \ldots
$$
and an asymptotic formula is known to hold \cite{McKay}
$$
\text{Vol}(\mathcal{B}_k) \sim \sqrt{2\pi} e^{1/3} \cdot \frac{k^{-(k-1)^2} e^{k^2}}{(2\pi)^{k}} \ , \ k \rightarrow \infty.
$$
In addition the asymptotic behavior of the Euler product $c_k$ is known (see \cite[Proposition]{ConreyGonek}), 
$$
\log c_k = -k^2 \log (2 e^{\gamma} \log k) + o(k^2)
$$
where $\gamma$ is the Euler--Mascheroni constant. 

We note that Conrey and Gamburd~\cite{congamb} compute the even integer moments $$\lim_{T \rightarrow \infty} \frac{1}{T}
\int_{0}^{T} \Bigg | \sum_{n \leq N} n^{-1/2 - it} \Bigg |^{2k} dt$$ on the $\tfrac 12$-line, and unsurprisingly the answer is extremely similar to Theorem \ref{thm:asymp} (in particular an Euler product and a volume related to the Birkhoff polytope again appear). Conrey and Gamburd discuss the connection between their result and the moments of certain truncated characteristic polynomials of random matrices. In general, it seems reasonable to say that the arithmetic factor $c_k$ reflects local counting modulo different primes in the moment computation, whereas the geometric factor $\gamma_k$ reflects global counting of tuples $n_1 , ... , n_k$ and $m_1 , ... , m_k$ subject to the truncation $n_i , m_i \leq N$.

We deduce Theorem \ref{thm:asymp} from a general result of La Bret\`{e}che \cite{breteche} on mean values
of multiplicative functions in several variables. 
 Theorem 3 has also been obtained independently by Granville and Soundararajan (unpublished), and also very recently (and independently)
by Heap and Lindqvist \cite{Heap}. 
Additionally, Theorem 3 sheds light on the conjectural behavior of
moments of the theta functions, 
$$
\theta(x, \chi) = \sum_{n \geq 1} \chi(n) e^{-\pi n^2 x / p}
$$
with $p \geq 3$ a prime, and $\chi$ an even Dirichlet character modulo $p$.  
The rapidly decaying factor $ e^{-\pi n^2 x / p}$ essentially restricts the sum to those $n$ less than about $\sqrt{p}$ (if $x=1$, say), and the average behavior of
$\chi(n)$ with $n \ll p^{1/2}$ is similar to that of a
Steinhaus random multiplicative function. Therefore Theorem 3 leads to the conjecture that
$$
\frac{1}{p} \sum_{\substack{\chi \mod p \\ \chi \text{ even}}} |\theta(1, \chi)|^{2k}
\sim C_k p^{k/2} (\log p)^{(k-1)^2} \;\;\;\;\; \text{as} \; p \rightarrow \infty .
$$
In unpublished recent work the same conjecture was stated by Marc Munsch on the basis of his lower bound for
moments of $\theta(1, \chi)$. Louboutin and Munsch~\cite{MunschLouboutin} prove the conjecture for $k = 1$ and $k = 2$.
%

Combining Theorem \ref{thm:lpmain} and Theorem \ref{thm:asymp} suggests the following ``counter-conjecture''
to Helson's claim (\ref{helson}).  
\begin{conjecture1} \label{theconjecture}
If $f(n)$ is a Steinhaus random multiplicative function, then we have as $N \rightarrow \infty$, 
$$
\mathbb{E} \Bigg | \sum_{n \leq N} f(n) \Bigg |^{2q} \sim \begin{cases}
C(q) N^q, & \text{ for } 0 \leq q \leq 1 \\
C(q) N^q (\log N)^{(q-1)^2},  & \text{ for } 1 \leq q .
\end{cases}
$$
\end{conjecture1}
Conjecture \ref{theconjecture} 
suggests a possible line of attack on the problem of showing that
for a positive proportion of even characters $\chi$ modulo $p$, 
we have $\theta(1, \chi) \neq 0$. This would be based on comparing the
first and second \textit{absolute} moments, i.e
$$
\sum_{\substack{\chi \mod p \\ \chi \text{ even}}} |\theta(1, \chi)| \ \ \text{and} \ \sum_{\substack{\chi \mod p \\ \chi \text{ even}}} |\theta(1,\chi)|^2.
$$

We emphasise that we do not have a lot of evidence towards Conjecture 1 when $q \notin \mathbb{N}$, and perhaps especially when $0 < q < 1$, 
and it is conceivable the behaviour could be more complicated. 
However this is the simplest possible conjecture respecting the 
information that we now have. In addition for $q > 1$ it perhaps seems unlikely that
the distribution of the tails of $\sum_{n < N} f(n)$ (in a large
deviation regime) fluctuates so significantly that it would affect
the exponent $(q-1)^2$ of the logarithm when $q$ goes from an
integer to a fractional exponent. We also note that if we could obtain
the order of magnitude for the $2q$-th moment suggested by the Conjecture \ref{theconjecture}
for $q=\tfrac 12$, then since we know it trivially for $q=1$ a simple
argument using H\"older's inequality (as in the proof of Theorem
\ref{thm:lpmain}, below) would establish the order of magnitude 
suggested by the Conjecture for all $0 \leq q \leq 1$.

Finally, following a question from the referee, 
we noticed that we can extend Theorem 3 to the Rademacher case.
We omit the simple cases of $k = 1, 2$ in the theorem below, since
both are different from the case $k \geq 3$. 
\begin{theorem} \label{thm:rademachermoments}
Let $f(n)$ be a Rademacher random multiplicative function. Then, for $k \geq 3$
an integer,
as $N
\rightarrow \infty$, 
$$
\mathbb{E} \Big ( \sum_{n \leq N} f(n) \Big )^{k} 
\sim C_k  \cdot N^{k/2} (\log N)^{k(k-3)/2}
$$
with $C_k > 0$ constant. 
\end{theorem}
Similarly as in Theorem \ref{thm:asymp} the constant $C_k$ splits into
an arithmetic and geometric factor. The interested reader should have
no trouble working out the details. Theorem \ref{thm:rademachermoments}
has also been obtained independently by Heap and Lindqvist \cite{Heap}.

At first glance it may seem strange that all the moments here (including the odd ones) are non-trivially large, but that is because in the Rademacher case there is no distinction between a term and its complex conjugate (and similarly if one calculated an expression like $\mathbb{E} \left| \sum_{n \leq N} f(n) \right|^{2k} \left(\sum_{n \leq N} f(n) \right)$ in the Steinhaus case, this would be non-trivially large provided $k \geq 1$). Note also that the moments are rather larger in the Rademacher case than the Steinhaus case, again because everything is real valued and so the terms exhibit less cancellation.

\textbf{Acknowledgments} We are grateful to the referee for a careful reading of the paper and for asking several questions which led to 
Theorem 4 and stronger results in Theorem 3. 

\section{Lower bounds for the first moment}
In this section we shall first prove the following result.
\begin{proph1}
Let $f(n)$ be a Rademacher random multiplicative function. There exist arbitrarily large values of $x$ for which
$$ \EE\left| \sum_{n \leq x} f(n) \right| \geq \frac{\sqrt{x}}{(\log\log x)^{3+o(1)}} . $$

The same is true if $f(n)$ is a Steinhaus random multiplicative function.
\end{proph1}

The above proposition is actually a fairly straightforward deduction from the work of Harper~\cite{harpergp}. However, it is a bit unsatisfactory because it only gives a lower bound along some special sequence of $x$ values. With more work we can correct this defect, as in the following theorem announced in the Introduction:
\begin{thmh1}
Let $f(n)$ be a Rademacher random multiplicative function. Then for {\em all} large $x$ we have
$$ \EE\left| \sum_{n \leq x} f(n) \right| \geq \frac{\sqrt{x}}{(\log\log x)^{3+o(1)}} . $$

The same is true if $f(n)$ is a Steinhaus random multiplicative function.
\end{thmh1}

The proof of Proposition 1 has two ingredients. The first is the observation, essentially due to Hal\'{a}sz~\cite{halasz}, that one can almost surely lower bound an average of $\left| \sum_{n \leq x} f(n) \right|$ in terms of the behaviour of $f(n)$ on primes only: more specifically, in the Rademacher case we almost surely have that, for any $y \geq 2$,
$$ \int_{1}^{\infty} \frac{\left| \sum_{n \leq z} f(n) \right|}{z^{3/2 + 1/\log y}} dz \gg \sup_{t \geq 1} \ \exp
\Big ( \sum_{p}\frac{f(p) \cos(t\log p)}{p^{1/2+1/\log y}} - \log t - \log\log(t+2)/2 \Big )  . $$
Here the implicit constant in the $\gg$ notation is absolute. The reader should note that the presence of the supremum over $t$ will be very significant here, since at any fixed $t$ the expected size of the right hand side would be too small to produce a useful result (about $\log^{1/4}y$, rather than about $\log y$ which is what we need).

The second ingredient is a strong lower bound for the expected size of the right hand side, which we deduce from the work of Harper~\cite{harpergp}. We quote the relevant statements from Harper's work as a lemma now, for ease of reference later.
\begin{lemh1}\textup{(See $\S 6.3$ of \cite{harpergp}.)}
If $(f(p))_{p \; \text{prime}}$ are independent Rademacher random variables, then with probability $1-o(1)$ as $x \rightarrow \infty$ we have
\begin{align*} \sup_{1 \leq t \leq 2(\log\log x)^{2}} \sum_{p} \frac{f(p) \cos(t\log p)}{p^{1/2+1/\log x}} \geq \log\log x - \log & \log\log x  \\  - & O((\log\log\log x)^{3/4}) . \end{align*}

If $(f(p))_{p \; \text{prime}}$ are independent Steinhaus random variables, then with probability $1-o(1)$ as $x \rightarrow \infty$ we have
\begin{align*} \sup_{1 \leq t \leq 2(\log\log x)^{2}} \sum_{p} \Big ( \frac{\Re(f(p)p^{-it})}{p^{1/2 + 1/\log x}} + \frac{1}{2} \frac{\Re(f(p)^{2}p^{-2it})}{p^{1 + 2/\log x}} \Big ) \geq \log & \log x -  \log\log\log x \\- & O((\log\log\log x)^{3/4}) . \end{align*}
\end{lemh1}
The first statement here is proved in the last paragraph in $\S 6.3$ of \cite{harpergp} (noting that the quantity $y$ there is $\log^{8}x$). The second statement can be proved by straightforward adaptation of that argument, the point being that the expectation and covariance structure of these random sums in the Steinhaus case are the same, up to negligible error terms, as in the Rademacher case, so the same arguments can be applied. (See the preprint \cite{harpertypicalmax} for an explicit treatment of some very similar Steinhaus random sums.) The argument in \cite{harpergp} is quite involved, but the basic aim is to show that, for the purpose of taking the supremum, the sums $\sup_{1 \leq t \leq 2(\log\log x)^{2}} \sum_{p} \frac{f(p) \cos(t\log p)}{p^{1/2+1/\log x}}$ behave somewhat independently at values of $t$ that are separated by $\gg 1/\log x$, so one has something like the supremum over $\log x$ independent samples.

To prove Theorem 1 we introduce a third ingredient, namely we show that $\EE\left| \sum_{n \leq x} f(n) \right|$ may itself be lower bounded in terms of an integral average of $\EE\left| \sum_{n \leq z} f(n) \right|$, as follows:
\begin{proph2}
Let $f(n)$ be a Rademacher random multiplicative function. For any large $x$ we have
$$ \EE\left|\sum_{n \leq x} f(n) \right| \gg \frac{\sqrt{x}}{\log x} \int_{1}^{\sqrt{x}} \left(\frac{\EE|\sum_{n \leq z} f(n)|}{\sqrt{z}} \right) \frac{dz}{z} . $$

The same is true if $f(n)$ is a Steinhaus random multiplicative function.
\end{proph2}
This uses the multiplicativity of $f(n)$ in an essential way (as does the proof of Proposition 1, of course).

Theorem 1 then follows quickly by combining Proposition 2 with the proof of Proposition 1.

\vspace{12pt}
As the reader will see, the proof of Proposition 2 is based on a ``physical space'' decomposition of the sum $\sum_{n \leq x} f(n)$, which is somewhat related to the martingale arguments of Harper~\cite{harperlimits}. This is unlike the other arguments above, which work by establishing a connection between the integral average of $\sum_{n \leq x} f(n)$ and its Dirichlet series $\sum_{n} f(n)/n^{s}$ (on the ``Fourier space'' side).

\subsection{Proof of Proposition 1}
The proof of Proposition 1 is slightly cleaner in the Rademacher case, because then $f(p)^{2} \equiv 1$ for all primes $p$. So we shall give the proof in that case first, and afterwards explain the small changes that arise in the Steinhaus case.

We know from work of Wintner~\cite{wintner} that almost surely $\sum_{n \leq x} f(n) = O_{\epsilon}(x^{1/2+\epsilon})$. Consequently, by partial summation the Dirichlet series $F(s) := \sum_{n} f(n)/n^{s}$ is almost surely convergent in the half plane $\Re(s) > 1/2$, and then by term by term integration it satisfies
$$ F(s) = s \int_{1}^{\infty} \frac{\sum_{n \leq z} f(n)}{z^{s+1}} dz , \;\;\;\;\; \Re(s) > 1/2 . $$
In particular, $F(s)$ is almost surely a holomorphic function on the half plane $\Re(s) > 1/2$.

On the other hand, since $f(n)$ is multiplicative we have for any $\Re(s) > 1$ that, in the Rademacher case,
\begin{align*}
F(s) = \prod_{p} \left(1+ \frac{f(p)}{p^{s}} \right) & = \exp \Big ( \sum_{p} \log\left(1+ \frac{f(p)}{p^{s}}\right) \Big )   \\
=  \exp & \Big ( \sum_{p} \frac{f(p)}{p^{s}} - \frac{1}{2} \sum_{p} \frac{f(p)^{2}}{p^{2s}} + \sum_{k \geq 3} \frac{(-1)^{k+1}}{k} \sum_{p} \frac{f(p)^{k}}{p^{ks}} \Big ) . 
\end{align*}
Therefore in the Rademacher case we have
$$ s \int_{1}^{\infty} \frac{\sum_{n \leq z} f(n)}{z^{s+1}} dz = \exp \Big ( \sum_{p} \frac{f(p)}{p^{s}} - \frac{1}{2} \sum_{p} \frac{f(p)^{2}}{p^{2s}} + \sum_{k \geq 3} \frac{(-1)^{k+1}}{k} \sum_{p} \frac{f(p)^{k}}{p^{ks}} \Big ) $$
at least when $\Re(s) > 1$, since both sides are equal to $F(s)$. But all the sums involving $p^{2s}$ and $p^{ks}$ are clearly absolutely convergent whenever $\Re(s) > 1/2$, and therefore define holomorphic functions there. In addition, for any fixed $s$ with $\Re(s) > 1/2$ the series $\sum_{p} \frac{f(p)}{p^{s}}$ is a sum of independent random variables, and Kolmogorov's Three Series Theorem implies it converges almost surely. Since a Dirichlet series is a holomorphic function strictly to the right of its abscissa of converge, we find that almost surely $\sum_{p} \frac{f(p)}{p^{s}}$ is a holomorphic function on the half plane $\Re(s) > 1/2$, and so almost surely we have, for all $\Re s > \tfrac 12$, 
$$ s \int_{1}^{\infty} \frac{\sum_{n \leq z} f(n)}{z^{s+1}} dz = \exp \Big ( \sum_{p} \frac{f(p)}{p^{s}} - \frac{1}{2} \sum_{p} \frac{f(p)^{2}}{p^{2s}} + \sum_{k \geq 3} \frac{(-1)^{k+1}}{k} \sum_{p} \frac{f(p)^{k}}{p^{ks}} \Big )  . $$
Next, if we write $s = \sigma + it$ and take absolute values on both sides then we find that, almost surely,
\begin{align*}
|s| \int_{1}^{\infty} \frac{\left|\sum_{n \leq z} f(n)\right|}{z^{\sigma+1}} dz & \geq \exp \Bigg ( \Re\Big(\sum_{p} \frac{f(p)}{p^{s}} - \frac{1}{2} \sum_{p} \frac{f(p)^{2}}{p^{2s}} + \sum_{k \geq 3} \frac{(-1)^{k+1}}{k} \sum_{p} \frac{f(p)^{k}}{p^{ks}} \Big ) \Bigg ) \\
 =  \exp \Big ( \sum_{p} & \frac{\Re(f(p)p^{-it})}{p^{\sigma}} - \frac{1}{2} \sum_{p} \frac{\Re(f(p)^{2}p^{-2it})}{p^{2\sigma}} + O(1) \Big ),  \;\;\;\;\; \forall \sigma > 1/2 . 
\end{align*}
If we take $\sigma = 1/2 + 1/\log y$ for a parameter $y \geq 2$, and we note that then $|s| \asymp |t|$ provided $t \geq 1$ (say), we have almost surely that for all $y \geq 2$, 
$$ \int_{1}^{\infty} \frac{\left|\sum_{n \leq z} f(n)\right|}{z^{3/2 + 1/\log y}} dz \gg \sup_{t \geq 1} \exp \Big ( \sum_{p} \frac{\Re(f(p)p^{-it})}{p^{1/2 + 1/\log y}} - \frac{1}{2} \sum_{p} \frac{\Re(f(p)^{2}p^{-2it})}{p^{1 + 2/\log y}} - \log t \Big ) . $$
In the Rademacher case the first sum over $p$ is $\sum_{p} \frac{f(p)\cos(t\log p)}{p^{1/2 + 1/\log y}}$, and (since $f(p)^{2}=1$) the second sum over $p$ is $\Re \sum_{p} \frac{1}{p^{1 + 2/\log y + 2it}} = \Re \log\zeta(1+2/\log y + 2it) + O(1)$, where $\zeta$ denotes the Riemann zeta function. Standard estimates (see e.g. Theorem 6.7 of Montgomery and Vaughan~\cite{mv}) imply that $|\log\zeta(1+2/\log y + 2it)| \leq \log\log(t+2) + O(1)$ for $t \geq 1$, so we have almost surely that for all $y \geq 2$, 
\begin{equation}\label{halineq}
\int_{1}^{\infty} \frac{\left|\sum_{n \leq z} f(n)\right|}{z^{3/2 + 1/\log y}} dz \gg \sup_{t \geq 1} \exp \Big ( \sum_{p} \frac{f(p)\cos(t\log p)}{p^{1/2 + 1/\log y}} - \log t - \log\log(t+2)/2 \Big ). 
\end{equation}

(The above argument and inequality \eqref{halineq} are essentially due to Hal\'{a}sz~\cite{halasz}, and are also related to the arguments of Wintner~\cite{wintner}. The only small difference is that Hal\'{a}sz restricted to $1 \leq t \leq 2$. See Appendix A of Harper~\cite{harpergp} for a presentation similar to the above.)  

\vspace{12pt}
Now to prove Proposition 1, note that for any large parameters $x$ and $x_0 < x_1$ we have
$$ \sup_{x_0 < z < x_1} \frac{\EE|\sum_{n \leq z} f(n)|}{\sqrt{z}} \geq \frac{1}{\log x} \int_{x_0}^{x_1} \frac{\EE\left|\sum_{n \leq z} f(n)\right|}{z^{3/2 + 1/\log x}} dz , $$
since $\int_{x_0}^{x_1} \frac{dz}{z^{1+1/\log x}} \leq \int_{1}^{\infty} \frac{dz}{z^{1+1/\log x}} = \log x$. Then by Cauchy--Schwarz we always have $\EE|\sum_{n \leq z} f(n)| \leq \sqrt{z}$, so
\begin{eqnarray}
\int_{x_0}^{x_1} \frac{\EE\left|\sum_{n \leq z} f(n)\right|}{z^{3/2 + 1/\log x}} dz & \geq & \int_{1}^{\infty} \frac{\EE\left|\sum_{n \leq z} f(n)\right|}{z^{3/2 + 1/\log x}} dz - \int_{1}^{x_0} \frac{dz}{z^{1+1/\log x}} - \int_{x_1}^{\infty} \frac{dz}{z^{1+1/\log x}} \nonumber \\
& \geq & \int_{1}^{\infty} \frac{\EE\left|\sum_{n \leq z} f(n)\right|}{z^{3/2 + 1/\log x}} dz - \log x_0 - \frac{\log x}{x_1^{1/\log x}} . \nonumber
\end{eqnarray}
In particular, if we choose $x_0 = e^{\sqrt{\log x}}$ and $x_1 = e^{(\log x) \log\log x}$, say, then we have
\begin{equation}\label{divineq}
\sup_{x_0 < z < x_1} \frac{\EE|\sum_{n \leq z} f(n)|}{\sqrt{z}} \geq \frac{1}{\log x} \int_{1}^{\infty} \frac{\EE\left|\sum_{n \leq z} f(n)\right|}{z^{3/2 + 1/\log x}} dz - \frac{2}{\sqrt{\log x}} .
\end{equation}

Finally, in the Rademacher case Lemma 1 implies that, with probability $1-o(1)$ as $x \rightarrow \infty$,
$$ \sup_{1 \leq t \leq 2(\log\log x)^{2}} \sum_{p} \frac{f(p) \cos(t\log p)}{p^{1/2+1/\log x}} \geq \log\log x - \log\log\log x - O((\log\log\log x)^{3/4}) . $$
This implies that with probability $1-o(1)$ one has
$$ \sup_{t \geq 1} \exp \Big ( \sum_{p} \frac{f(p)\cos(t\log p)}{p^{1/2 + 1/\log x}} - \log t - \log\log(t+2)/2 \Big ) \geq \frac{\log x}{(\log\log x)^{3+o(1)}} , $$
and then by the Hal\'{a}sz type lower bound inequality \eqref{halineq} we deduce
\begin{equation}\label{intexbound}
\int_{1}^{\infty} \frac{\EE\left|\sum_{n \leq z} f(n)\right|}{z^{3/2 + 1/\log x}} dz \geq \frac{\log x}{(\log\log x)^{3+o(1)}} .
\end{equation}
Proposition 1 follows in the Rademacher case by combining this with \eqref{divineq}.

\vspace{12pt}
In the Steinhaus case the initial argument of Wintner~\cite{wintner} still works, so the first change that is needed in the preceding argument comes in the expression for the Euler product $F(s)$, which for $\Re(s) > 1$ is now
\begin{eqnarray}
F(s) = \prod_{p} \left(1+ \sum_{j=1}^{\infty} \frac{f(p)^{j}}{p^{js}}\right) & = & \exp \Big ( - \sum_{p} \log\left(1- \frac{f(p)}{p^{s}}\right) \Big ) \nonumber \\
& = & \exp \Big ( \sum_{p} \frac{f(p)}{p^{s}} + \frac{1}{2} \sum_{p} \frac{f(p)^{2}}{p^{2s}} + \sum_{k \geq 3} \frac{1}{k} \sum_{p} \frac{f(p)^{k}}{p^{ks}} \Big ) . \nonumber
\end{eqnarray}
Notice this is the same as we had in the Rademacher case, except now there are no alternating minus signs in the final exponential. The argument using the Three Series Theorem, etc. then continues as in the Rademacher case to yield that, almost surely,
$$ s \int_{1}^{\infty} \frac{\sum_{n \leq z} f(n)}{z^{s+1}} dz = \exp \Big ( \sum_{p} \frac{f(p)}{p^{s}} + \frac{1}{2} \sum_{p} \frac{f(p)^{2}}{p^{2s}} + \sum_{k \geq 3} \frac{1}{k} \sum_{p} \frac{f(p)^{k}}{p^{ks}} \Big ) \;\;\;\;\; \forall \; \Re(s) > 1/2 . $$
Putting $s=1/2 + 1/\log y + it$ and taking absolute values on both sides, we deduce that almost surely,
\begin{equation}\label{halineq2}
\int_{1}^{\infty} \frac{\left|\sum_{n \leq z} f(n)\right|}{z^{3/2 + 1/\log y}} dz \gg \sup_{t \geq 1} \ \exp \Big ( \sum_{p} ( \frac{\Re(f(p)p^{-it})}{p^{1/2 + 1/\log y}} + \frac{1}{2} \frac{\Re(f(p)^{2}p^{-2it})}{p^{1 + 2/\log y}}) - \log t \Big ) \ ,  \forall y \geq 2 .
\end{equation}

Since we don't now have $f(p)^{2} \equiv 1$, we cannot remove the contribution of the prime squares using estimates for the zeta function. However, by the Steinhaus case of Lemma 1 we still have that, with probability $1-o(1)$ as $x \rightarrow \infty$,
\begin{align*} \sup_{1 \leq t \leq 2(\log\log x)^{2}} \sum_{p} \Big ( \frac{\Re(f(p)p^{-it})}{p^{1/2 + 1/\log x}} + \frac{1}{2} \frac{\Re(f(p)^{2}p^{-2it})}{p^{1 + 2/\log x}} \Big ) \geq \log & \log x - \log\log\log x \\ - & O((\log\log\log x)^{3/4}) , \end{align*}
and therefore with probability $1-o(1)$ we have
$$ \sup_{t \geq 1} \ \exp \Big ( \sum_{p} \Big ( \frac{\Re(f(p)p^{-it})}{p^{1/2 + 1/\log y}} + \frac{1}{2} \frac{\Re(f(p)^{2}p^{-2it})}{p^{1 + 2/\log y}} \Big ) - \log t \Big ) \geq \frac{\log x}{(\log\log x)^{3+o(1)}} . $$
Combining this estimate with \eqref{halineq2} and \eqref{divineq} then proves Proposition 1 in the Steinhaus case.

\subsection{Proofs of Theorem 1 and Proposition 2}

\begin{proof}[Proof of Theorem 1, assuming Proposition 2]
In view of Proposition 2, it will suffice to prove that for all large $x$ we have
$$ \int_{1}^{\sqrt{x}} \left(\frac{\EE|\sum_{n \leq z} f(n)|}{\sqrt{z}} \right) \frac{dz}{z} \geq \frac{\log x}{(\log\log x)^{3+o(1)}} . $$
However, for any large parameter $y$ we have
\begin{align*} \int_{1}^{\sqrt{x}} \left(\frac{\EE|\sum_{n \leq z} f(n)|}{\sqrt{z}} \right) \frac{dz}{z} & \geq \int_{1}^{\sqrt{x}} \frac{\EE|\sum_{n \leq z} f(n)|}{z^{3/2+1/\log y}} dz \\ & \geq  \frac{\log y}{(\log\log y)^{3+o(1)}} -  \int_{\sqrt{x}}^{\infty} \frac{\EE|\sum_{n \leq z} f(n)|}{z^{3/2+1/\log y}} dz , 
\end{align*}
in view of the lower bound $\int_{1}^{\infty} \frac{\EE|\sum_{n \leq z} f(n)|}{z^{3/2+1/\log y}} dz \geq \frac{\log y}{(\log\log y)^{3+o(1)}}$ obtained in \eqref{intexbound}. By Cauchy--Schwarz we always have $\EE|\sum_{n \leq z} f(n)| \leq \sqrt{z}$, so the subtracted term here is at most
$$ \int_{\sqrt{x}}^{\infty} \frac{dz}{z^{1+1/\log y}} = \frac{\log y}{(\sqrt{x})^{1/\log y}} . $$
If we choose $\log y$ somewhat smaller than $\log x$, say $\log y = (\log x)/(100\log\log\log x)$, we deduce that
$$ \int_{1}^{\sqrt{x}} \left(\frac{\EE|\sum_{n \leq z} f(n)|}{\sqrt{z}} \right) \frac{dz}{z} \geq \frac{\log x}{(\log\log x)^{3+o(1)}} - \frac{\log x}{(\log\log x)^{50}} = \frac{\log x}{(\log\log x)^{3+o(1)}} , $$
as required.
\end{proof}

\begin{proof}[Proof of Proposition 2]
The first part of the proof again differs slightly depending on whether we are in the Rademacher or the Steinhaus case. We will first work in the Rademacher case and then explain the small changes needed in the other situation.

Let $A_t := \sum_{n \leq t} f(n)$. If we let $P(n)$ denote the largest prime factor of $n$, we have
$$ \sum_{n \leq x} f(n) = \sum_{p \leq x} \sum_{n \leq x, P(n)=p} f(n) = \sum_{p \leq x} f(p) \sum_{m \leq x/p, P(m) < p} f(m) , $$
since $f$ is multiplicative. Here the inequality $P(m) < p$ in the final sum is strict because $f$ is supported on squarefree numbers. Notice here that if $p > \sqrt{x}$ then $x/p < \sqrt{x} < p$, so we automatically have $P(m) < p$ in the inner sums over $m$. Thus we can rewrite things slightly as
\begin{align*} \sum_{n \leq x} f(n) & = \sum_{\sqrt{x} < p \leq x} f(p) \sum_{m \leq x/p} f(m) + \sum_{p \leq \sqrt{x}} f(p) \sum_{m \leq x/p, P(m) < p} f(m) \\ & =: \sum_{\sqrt{x} < p \leq x} f(p) A_{x/p} + B_x , \end{align*}
say. Notice also that the random variables $A_{x/p}$ and $B_x$ are independent of the $f(p)$ for $\sqrt{x} < p \leq x$.

We shall introduce a penultimate piece of notation, by defining the random variable
$$ C_{x} := \sum_{\sqrt{x} < p \leq x} f(p) A_{x/p} . $$
Finally, let $\epsilon$ be a Rademacher random variable that is independent of everything else.

Now since the $(f(p))_{\sqrt{x} < p \leq x}$ are symmetric random variables independent of $B_x$ and the $A_{x/p}$, it follows that 
$$ \sum_{n \leq x} f(n) = \sum_{\sqrt{x} < p \leq x} f(p) A_{x/p} + B_x \stackrel{d}{=} \epsilon \sum_{\sqrt{x} < p \leq x} f(p) A_{x/p} + B_x , $$
where $\stackrel{d}{=}$ denotes equality in distribution. Then if we {\em condition on the values of $B_x , C_x$}, we find the conditional expectation
$$ \EE\left( \left|\epsilon \sum_{\sqrt{x} < p \leq x} f(p) A_{x/p} + B_x \right| \Bigg| B_x , C_x \right) = (1/2) |C_x + B_x| + (1/2) |-C_x + B_x| \geq |C_x| , $$
by the triangle inequality. Now if we average over values of $B_x, C_x$, and use the Tower Property of conditional expectations (the fact that the expectation of a conditional expectation is the unconditional expectation), we obtain
$$ \EE\left|\sum_{n \leq x} f(n) \right| = \EE\left|\epsilon \sum_{\sqrt{x} < p \leq x} f(p) A_{x/p} + B_x \right| \geq \EE|C_x| . $$

On recalling the definitions of $C_x$ and $A_{x/p}$, we see we have proved the following:
\begin{lemh2}
For all large $x$ we have
$$ \EE\left|\sum_{n \leq x} f(n)\right| \geq \EE\left| \sum_{\sqrt{x} < p \leq x} f(p) \sum_{m \leq x/p} f(m) \right| . $$
\end{lemh2}

(In the Steinhaus case one has a weak inequality $P(m) \leq p$ in the definition of $B_x$, since $f$ is totally multiplicative, but this makes no difference to the argument just given. Instead of choosing $\epsilon$ to be a Rademacher random variable one can choose $\epsilon$ to be uniformly distributed on the unit circle, and then one obtains exactly the same conclusion in Lemma 2.)

\vspace{12pt}
Since the $f(p)$ are Rademacher or Steinhaus random variables independent of the ``coefficients'' $\sum_{m \leq x/p} f(m) = A_{x/p}$, an application of Khintchine's inequality (see e.g. Gut's textbook~\cite{gut}) yields that
$$ \EE\left| \sum_{\sqrt{x} < p \leq x} f(p) \sum_{m \leq x/p} f(m) \right| \gg \EE\sqrt{\sum_{\sqrt{x} < p \leq x} \left|\sum_{m \leq x/p} f(m)\right|^{2}} . $$
It would be nice if we could find a way to exploit this (sharp) bound with the squares still in place on the inside, but to prove Proposition 2 we shall trade them away in order to remove the intractable squareroot. Thus by the Cauchy--Schwarz inequality and the fact that $\sum_{\sqrt{x} < p \leq x} 1/p = \log 2 + o(1)$ we have
\begin{align*} \sum_{\sqrt{x} < p \leq x} \sqrt{\frac{1}{p}} \left|\sum_{m \leq x/p} f(m)\right| & \leq \sqrt{\sum_{\sqrt{x} < p \leq x} \frac{1}{p}} \sqrt{\sum_{\sqrt{x} < p \leq x} \left|\sum_{m \leq x/p} f(m)\right|^{2}} \\ & \ll \sqrt{\sum_{\sqrt{x} < p \leq x} \left|\sum_{m \leq x/p} f(m)\right|^{2}} . 
\end{align*}
Combining this with the above, we deduce:
\begin{lemh3}
For all large $x$ we have
$$ \EE\left|\sum_{n \leq x} f(n)\right| \gg \sum_{\sqrt{x} < p \leq x} \frac{1}{\sqrt{p}} \EE\left| \sum_{m \leq x/p} f(m) \right| \geq \frac{1}{\log x} \sum_{\sqrt{x} < p \leq x} \frac{\log p}{\sqrt{p}} \cdot \EE\left| \sum_{m \leq x/p} f(m) \right| . $$
\end{lemh3}

\vspace{12pt}
We have now almost finished the proof of Proposition 2. If we have two primes $z \leq p \leq p' \leq z+ z/\log^{1000}x$ for some $\sqrt{x} < z \leq x$ then
\begin{align*} \left| \EE\left| \sum_{m \leq x/p} f(m) \right| - \EE\left| \sum_{m \leq x/p'} f(m) \right| \right| & \leq \EE\left| \sum_{x/p' < m \leq x/p} f(m) \right| \\ & \ll \sqrt{x(\frac{1}{p} - \frac{1}{p'}) + 1} \ll \sqrt{\frac{x}{p \log^{1000}x}} + 1 , \end{align*}
by the Cauchy--Schwarz inequality and orthogonality of the $f(m)$. And we see
$$ \frac{1}{\log x} \sum_{\sqrt{x} < p \leq x} \frac{\log p}{\sqrt{p}} \left( \sqrt{\frac{x}{p \log^{1000}x}} + 1 \right) \ll \frac{\sqrt{x}}{\log^{500}x} + \frac{1}{\log x} \sum_{\sqrt{x} < p \leq x} \frac{\log p}{\sqrt{p}} \ll \frac{\sqrt{x}}{\log x} , $$
which will make a negligible contribution in Proposition 2, so in Lemma 3 we may replace each term $\EE\left| \sum_{m \leq x/p} f(m) \right|$ by an averaged version
$$ \frac{\log^{1000}x}{p} \int_{p}^{p(1+1/\log^{1000}x)} \EE\left| \sum_{m \leq x/t} f(m) \right| dt . $$
Since we know that primes are well distributed in intervals of relative length $1 + 1/\log^{1000}x$ (with density 1 when weighted by $\log p$) we can rewrite Lemma 3 as
\begin{eqnarray}
\EE\left|\sum_{n \leq x} f(n)\right| & \gg & \frac{1}{\log x} \sum_{\sqrt{x} < p \leq x} \log p \frac{\log^{1000}x}{p} \int_{p}^{p(1+1/\log^{1000}x)} \EE\left| \sum_{m \leq x/t} f(m) \right| \frac{dt}{\sqrt{t}} \nonumber \\
& \gg & \frac{1}{\log x} \int_{\sqrt{x}}^{x} \EE\left| \sum_{m \leq x/t} f(m) \right| \frac{dt}{\sqrt{t}} . \nonumber
\end{eqnarray}
Proposition 2 now follows by making the substitution $z = x/t$ in the integral.
\end{proof}

\section{Lower bounds for small moments - Proof of Theorem 2}
The proof is a very simple argument using the Cauchy--Schwarz inequality and H\"{o}lder's inequality.

Indeed, for any $0 \leq q \leq 1$ we have
\begin{align*}
\mathbb{E} \Big | \sum_{n \leq N} f(n) \Big |
& \leq \mathbb{E} \Big [ \Big | \sum_{n \leq N} f(n) \Big |^{2q} \Big ]^{1/2}
\cdot \mathbb{E} \Big [ \Big | \sum_{n \leq N} f(n) \Big |^{2 - 2q} \Big ]^{1/2} \\
& \leq \mathbb{E} \Big [ \Big | \sum_{n \leq N} f(n) \Big |^{2q} \Big ]^{1/2}
\cdot \mathbb{E} \Big [ \Big | \sum_{n \leq N} f(n) \Big |^{2} \Big ]^{(1 - q)/2}.
\end{align*}
Since $\EE|\sum_{n \leq N} f(n)|^{2} \leq N$ and $\EE|\sum_{n \leq N} f(n)| \geq \sqrt{N}/(\log\log N)^{3+o(1)}$, by re-arranging we obtain the lower bound
$$
\mathbb{E} \Big [ \Big | \sum_{n \leq N} f(n) \Big |^{2q} \Big ] \geq N^{q} (\log\log N)^{-6 + o(1)}.
$$

\section{Asymptotics for even moments - Proof of Theorem 3}

Note that
\begin{align} \nonumber
\mathbb{E} \Big | \sum_{n \leq X} f(n) \Big |^{2k} & = \sum_{\substack{n_1, \ldots, n_k \leq X \\ m_1, \ldots, m_k \leq X}}
\mathbb{E} [ f(n_1) \ldots f(n_k) \overline{f(m_1) \ldots f(m_k)} ] \\
& \label{equation1} 
= \sum_{\substack{n_1, \ldots, n_k \leq X \\ m_1, \ldots, m_k \leq X \\ n_1 \ldots n_k = m_1 \ldots m_k}} 1 .
\end{align}
Now 
$$
g(n_1, \ldots, n_k, m_1, \ldots, m_k) = \mathbf{1}_{n_1 \ldots n_k = m_1 \ldots m_k}
$$
is a multiplicative function of several variables\footnote{In other words $$g(n_1, \ldots, n_k, m_1, \ldots, m_k) g(u_1, \ldots, u_k, v_1, \ldots, v_k) = g(n_1 u_1, \ldots, n_k u_k, m_1 v_1, \ldots, m_k v_k)$$ for any natural numbers $n_i , m_i$ and $u_i , v_i$ whose least common multiples are coprime.} and our problem reduces to understanding
the mean value of
$$
\sum_{\substack{n_1, \ldots, n_k \leq X \\ m_1, \ldots, m_k \leq X}} g(n_1, \ldots, n_k, m_1, \ldots, m_k) .
$$
We notice that the associated multiple Dirichlet series
$$
\sum_{\substack{n_1, \ldots, n_k \\ m_1, \ldots , m_k}} \frac{g(n_1, \ldots, n_k, m_1, \ldots, m_k)}{n_1^{s_1} \ldots n_k^{s_k}
m_1^{w_1} \ldots m_k^{w_k}} = \sum_{n} \sum_{n_1 n_2 \ldots n_k = n} \frac{1}{n_1^{s_1} \ldots n_k^{s_k}} \sum_{m_1 m_2 \ldots m_k = n} \frac{1}{m_1^{w_1} \ldots m_k^{w_k}}
$$
is absolutely convergent for $\Re s_i, \Re w_i > \tfrac 12$ and moreover it factors as
$$
H(s_1, \ldots, s_k, w_1, \ldots, w_k) \prod_{i=1}^{k} \prod_{j=1}^{k} \zeta(s_i + w_j) 
$$
with $H(s_1, \ldots, s_k, w_1, \ldots, w_k)$ absolutely convergent in the region $\Re s_i, \Re w_i > \tfrac 14$. In addition
a direct check shows that $$
H(\tfrac 12,\ldots, \tfrac 12) = \prod_{p} \Big ( 1 - \frac{1}{p}
\Big )^{k^2} \cdot \Big ( 1 + \frac{k^2}{p} + \sum_{\alpha \geq 2}
\frac{\binom{a + k - 1}{k - 1}^2}{p^{\alpha}} \Big ) > 0.$$
Therefore the main result of La Bret\`{e}che \cite{breteche}
is applicable with the $k^2$ linear forms $\ell^{(i,j)}(s_1, \ldots, s_k, w_1, \ldots, w_k) := s_i + w_j$ with $1 \leq i,j \leq k$. We note that
the rank of the collection of linear forms $\ell^{(i,j)}$ (inside the space of all $\C$-linear forms on $\C^{2k}$) is $2k -1$. Therefore it follows from La Bret\`{e}che's result that (\ref{equation1})
is equal to
$$
(1 + o(1)) C_k X^k (\log X)^{k^2 - (2k -1)} .
$$

Using Th\'eor\`eme 2 in La Bret\'{e}che's work allows us to recover the
precise value of $C_k$. Indeed, according to Th\'eor\`eme 2 in \cite{breteche}
we get that (\ref{equation1}) is equal to
$$
(1 + o(1)) H(\tfrac 12, \ldots, \tfrac 12) \text{Vol}(A_k(X))
$$ 
where
$A_k(X)$ is a subset of $[1,\infty)^{k^2}$ corresponding to tuples $(a_{i,j}) \in [1,\infty)^{k^2}$ 
with $1 \leq i , j \leq k$ such that
\begin{align*}
& \text{for each } j \leq k: \prod_{1 \leq i \leq k} a_{i,j} \leq X \\
\text{ and } & \text{for each } i \leq k: \prod_{1 \leq j \leq k} a_{i,j} \leq X
\end{align*}
Therefore it remains to understand the asymptotic behavior of
$$
\text{Vol}(A_k(X))
$$
as $X \rightarrow \infty$. Surprisingly, this is somewhat involved, and the rest of the proof is devoted to that. 
\begin{proph3}
Let $k \geq 2$ be fixed. Then,
$$
\text{Vol}(A_k(X)) \sim \binom{2k-2}{k-1} k^{-(k-1)} \cdot \text{Vol}(\mathcal{B}_k)
\cdot X^k \cdot (\log X)^{(k-1)^2}
$$ 
where $\text{Vol}(\mathcal{B}_k)$ corresponds to the $(k-1)^2$ dimensional volume
of the Birkhoff polytope $\mathcal{B}_k \subset \mathbb{R}^{k^2}$. 
\end{proph3} 
The proof of the Proposition depends on the following Lemma.
\begin{lemh4} \label{lem:mainlemma}
Let $n \geq 1$ be fixed. Then as $X \rightarrow \infty$ we have
$$ \iint_{\substack{0 \leq x_1, ..., x_n \leq \log X \\ 0 \leq y_1, \ldots, y_n \leq \log X}}
  \exp \Big ( \min (x_1 + ... + x_n, y_1 + ... + y_n )\Big ) dx_1 ... dy_n \sim \binom{2n}{n} X^{n} . $$
\end{lemh4}
\begin{proof}
Making the substitutions $v_i = \log X - x_i$ and $w_i = \log X - y_i$ in Lemma \ref{lem:mainlemma}, we see the integral there is the same as
$$ X^{n} \iint_{\substack{0 \leq v_1, ..., v_n \leq \log X \\ 0 \leq w_1, \ldots, w_n \leq \log X}} \exp \Big (- \max (v_1 + ... + v_n, w_1 + ... + w_n )\Big ) dv_1 ... dw_n . $$
Here we can extend all the ranges of integration up to positive infinity, at the cost of a multiplicative error term $1+o(1)$. Then by symmetry 
\begin{eqnarray}
&& \iint_{\substack{0 \leq v_1, ..., v_n \\ 0 \leq w_1, \ldots, w_n}} \exp\Big (- \max(v_1 + ... + v_n, w_1 + ... + w_n) \Big ) dv_1 ... dw_n \nonumber \\
& = & 2 \int_{0 \leq v_1, ..., v_n} \exp\Big (- (v_1 + ... + v_n) \Big ) \int_{w_1 + ... + w_n \leq v_1 + ... + v_n}  dv_1 ... dw_n , \nonumber
\end{eqnarray}
and making the further substitution $v= v_1 + ... + v_n$ in the integral, we see the above is
\begin{eqnarray}
& = & 2 \int_{0}^{\infty} e^{-v} \left(\int_{v_1 + ... + v_{n-1} \leq v}  dv_1 ... dv_{n-1} \right) \left(\int_{w_1 + ... + w_n \leq v}  dw_1 ... dw_n \right) dv \nonumber \\
& = & 2 \int_{0}^{\infty} e^{-v} v^{2n-1} \left(\int_{v_1 + ... + v_{n-1} \leq 1} dv_1 ... dv_{n-1} \right) \left(\int_{w_1 + ... + w_n \leq 1} dw_1 ... dw_n \right) dv . \nonumber
\end{eqnarray}
Here the two integrals in brackets are simply the volume of the standard $n-1$ simplex and the standard $n$ simplex, which are well known to be $1/(n-1)!$ and $1/n!$ respectively. Threfore the above integral is equal to
$$
\frac{2}{(n-1)! n!} \int_{0}^{\infty} e^{-v} v^{2n-1} dv = 2 \binom{2n-1}{n} = \binom{2n}{n}.
$$
We conclude that the integral in the statement of Lemma \ref{lem:mainlemma} is equal to
(as $X \rightarrow \infty$), 
$$
(1 + o(1)) \binom{2n}{n} X^n
$$
as claimed. 
\end{proof} 

We are now ready to prove the Proposition, and thus finish the proof of Theorem 3.
\begin{proof}[of Proposition 3]
Notice first that, if we set $u_{i,j} = \log a_{i,j}$, and if we write
$$ c_j = \sum_{1 \leq i \leq k} u_{i,j}  \text{ and }  r_i = \sum_{1 \leq j \leq k} u_{i,j} $$
for all $i,j \leq k$, then we find
$$ \text{Vol}(A_k(X)) = \int_{(u_{i,j})_{1 \leq i,j \leq k} \subseteq [0,\infty)^{k^{2}} : c_j , r_i \leq \log X \; \forall i,j \leq k} \exp\left(\sum_{i, j \leq k} u_{i,j}\right) du_{1,1} ... du_{k,k} . $$
To prove the proposition we shall obtain upper and lower bounds for the integral on the right that are asymptotically equal.

For convenience of writing, we start by introducing a little more notation. Let $ S_{k-1} := \sum_{i,j \leq k-1} u_{i,j}$. Let also $\mathcal{U}_{k, \varepsilon}(X)$ be the set of $u_{i,j}$ with $i,j \leq k-1$
for which
$$
\sum_{i \leq k-1} u_{i,j} \leq \log X \text{ and } \sum_{j \leq k-1} u_{i,j} \leq \log X \text{ and } \sum_{i,j \leq k-1} u_{i,j} > (k-2- \varepsilon) \log X.
$$
Considering the vector $\mathbf{u}$ of $u_{i,j}$ with $i,j \leq k-1$ as fixed, let $\mathcal{T}_{C, k}(\mathbf{u}, X)$ be the set of those $u_{k,i}$ with $i \leq k-1$
for which
$$
c_j \leq \log X \text{ for all } j \leq k-1
$$
Finally, again consider the $u_{i,j}$ with $i,j \leq k-1$ as fixed
let $\mathcal{T}_{R, k}(\mathbf{u}, X)$ be the set of those $u_{j,k}$ with $j \leq k-1$
for which
$$
r_i \leq \log X \text{ for all } i \leq k-1.
$$ 

We set $\epsilon = 1/\sqrt{\log X}$, say.
First seeking an upper bound, we note that if we have $S_{k-1} \leq (k-2-\epsilon)\log X$ then $S_k \leq (k-\epsilon)\log X$, and therefore the part of the integral where $S_{k-1} \leq (k-2-\epsilon)\log X$ contributes at most
$$ X^{k-\epsilon} \cdot \int_{(u_{i,j})_{1 \leq i,j \leq k} \subseteq [0,\infty)^{k^{2}} : c_j , r_i \leq \log X \; \forall i,j \leq k} 1 du_{1,1} ... du_{k,k} \leq X^{k-\epsilon} \log^{k^{2}}X . $$
This is asymptotically negligible (for any fixed $k$) by our choice of $\epsilon$. 
Meanwhile, the part of the integral where $S_{k-1} > (k-2-\epsilon)\log X$ is equal to
\begin{align} \label{equ:integral} & \int_{\mathcal{U}_{k, \varepsilon}(X)} \exp ( S_{k-1} )  \int_{\mathcal{T}_{C,k}(\mathbf{u}, X)} \exp \Big ( u_{k,1} + ... + u_{k,k-1} \Big ) \times  \\ \nonumber \times  & \int_{\mathcal{T}_{R,k}(\mathbf{u}, X)} \exp \Big ( u_{1,k} + ... + u_{k-1,k} \Big ) \int_{u_{k,k} : c_k , r_k \leq \log X} \exp ( u_{k,k} ) \ du_{1,1} ... du_{k,k} . 
\end{align}
Here the innermost integral is over those  $$0 \leq u_{k,k} \leq \log X - \max(u_{k,1} + ... + u_{k,k-1}, u_{1,k} + ... + u_{k-1,k}),$$ assuming the upper range of integration is at least zero. Therefore 
the innermost integral is certainly bounded above (extending the lower limit to negative infinity, and then performing the integration) by
$$ X \exp \Big (-\max(u_{k,1} + ... + u_{k,k-1}, u_{1,k} + ... + u_{k-1,k}) \Big ) . $$
Substituting this in, it follows that (\ref{equ:integral}) is less than
\begin{equation} \label{equ:integral2} X \int_{\mathcal{U}_{k,\varepsilon}(X)} \int_{\mathcal{T}_{C,k}(\mathbf{u}, X)} \int_{\mathcal{T}_{R,k}(\mathbf{u},X)} \exp\Big ( \min(\sum_{1 \leq j \leq k-1} c_j , \sum_{1 \leq i \leq k-1} r_i)\Big ) \prod_{(i,j) \neq (k,k)} d u_{i,j} . \end{equation}
At this point we change variables, letting $r_1,\ldots, r_{k-1}$ and $c_{1}, \ldots, c_{k-1}$ run through the interval
$[0,\log X]$
so that $u_{i,k} = r_i - \sum_{1 \leq j \leq k-1} u_{i,j}$ and $u_{k,j} = c_j - \sum_{1 \leq i \leq k-1} u_{i,j}$. 
Since $u_{i,k} \geq 0$ and $u_{k, j} \geq 0$ this change of variable implies the
additional condition that for all $i,j \leq k-1$, 
\begin{equation} \label{equ:cond1}
\sum_{j \leq k-1} u_{i,j} \leq r_i \text{ and } \sum_{i \leq k-1} 
u_{i,j} \leq c_j
\end{equation}
The Jacobian of this linear change of variable is equal to $1$ since the linear transformation taking the $(u_{i,j})$ with $(i,j) \neq (k,k)$ into $(r_\ell, c_\ell, u_{i,j})$ with $i,j,\ell \leq k-1$ is 
upper triangular with only $1$'s on the diagonal. 

Given $\mathbf{r} = (r_1, \ldots, r_{k-1})$ and $\mathbf{c} = (c_1, \ldots, c_{k-1})$ we let 
$\widetilde{\mathcal{U}}_{k,\varepsilon}(\mathbf{r}, \mathbf{c}, X)$ be the set of $u_{i,j}$ with $i, j \leq k-1$ satisfying the conditions (\ref{equ:cond1}) and the standing condition that
\begin{equation} \label{equ:cond2}
\sum_{i,j \leq k-1} u_{i,j} \geq (k-2 - \varepsilon) \log X ,
\end{equation}
and we let $\widetilde{\mathcal{T}}_k(X)$ be the set of
$0 \leq r_1 , \ldots, r_{k-1} \leq \log X$ and $0 \leq c_1 , \ldots, c_{k-1} \leq \log X$. 
Then (\ref{equ:integral2}) can be re-written as
\begin{equation} \label{equ:integral3}
X \int_{\widetilde{\mathcal{T}}_k(X)} 
\exp \Big ( \min ( \sum_{1 \leq j \leq k-1} c_j , \sum_{1 \leq i \leq k-1}
r_i ) \Big ) \text{Vol} \Big ( \widetilde{\mathcal{U}}_{k,\varepsilon}(\mathbf{r}, \mathbf{c}, X) \Big )
\prod_{i \leq k-1} d c_i \ d r_i
\end{equation}
Since $r_i, c_i \leq \log X$ for all $i \leq k-1$, we have
\begin{align*}
\text{Vol} & (\widetilde{\mathcal{U}}_{k,\varepsilon}(\mathbf{r}, \mathbf{c}, X))
\leq
\text{Vol}(\widetilde{\mathcal{U}}_{k,\varepsilon}(\mathbf{\log X, \log X, X})) \\
& = (\log X)^{(k-1)^2} \cdot \text{Vol}(\widetilde{\mathcal{U}}_{k,\varepsilon}(\mathbf{1}, \mathbf{1}, e)) \sim (\log X)^{(k-1)^2} \cdot \text{Vol}(\widetilde{\mathcal{U}}_{k,0}(\textbf{1}, \textbf{1}, e))
\end{align*}
as $X \rightarrow \infty$, where $\mathbf{\log X} := (\log X, \ldots, \log X)$
and $\mathbf{1} := (1, \ldots, 1)$, and where we recall for the final asymptotic that $\epsilon = 1/\sqrt{\log X}$. As already mentioned in the introduction
$\text{Vol}(\widetilde{\mathcal{U}}_{k,0}(\mathbf{1},\mathbf{1},e)) = k^{-(k-1)} \text{Vol}(\mathcal{B}_k)$ where $\mathcal{B}_k$ is the Birkhoff polytope. 
It follows that (\ref{equ:integral2}) is
\begin{align*}
\leq (1 + o(1)) X (\log X)^{(k-1)^2} & \cdot k^{-(k-1)} \text{Vol} ( \mathcal{B}_k)
 \\ & \times \int_{\widetilde{\mathcal{T}}_k(X)} \exp \Big ( \min( \sum_{1 \leq j \leq k-1}
c_j, \sum_{1 \leq i \leq k-1} r_i ) \Big ) \prod_{i \leq k-1} d c_i d r_i
\end{align*}
and by Lemma \ref{lem:mainlemma} this is less than or equal to to
$$
(1 + o(1)) X (\log X)^{(k-1)^2} \cdot k^{-(k-1)} \text{Vol} ( \mathcal{B}_k)
\cdot X^{k-1} \binom{2k-2}{k-1} 
$$
thus finishing the proof of the upper bound. 

For the lower bound we restrict attention, as we may (due to positivity), 
to the part of the integral  where 
$S_{k-1} > (k-2 + \varepsilon) \log X$ and each $r_i$, $c_i$
is $\geq (1 - \varepsilon) \log X$. The point of the former condition is that
if it is satisfied then
$$
u_{1,k} + u_{2,k} + \ldots + u_{k-1, k} \leq (k-1)\log X - S_{k-1} \leq
(1 - \varepsilon ) \log X
$$
and similarly $u_{k,1} + u_{k,2} + \ldots + u_{k,k-1} \leq (1 - \varepsilon) \log X$, and therefore
$$
\log X - \max \Big ( u_{k,1} + \ldots + u_{k,k-1}, u_{1,k} + \ldots + u_{k-1,k} \Big ) > \varepsilon \log X = \sqrt{\log X} \rightarrow \infty.
$$
Therefore arguing as above the innermost integral over $u_{k,k}$ in (\ref{equ:integral}) contributes
$$
(1 + o(1)) X \exp \Big ( - \max(u_{k,1} + \ldots + u_{k,k-1}, u_{1,k} + \ldots + u_{k-1,k} ) \Big ) .
$$
Proceeding as before we thus arrive to (\ref{equ:integral3}) but with the additional condition that $(1 - \varepsilon) \log X < r_i, c_i < \log X$ (and with the condition that $\sum_{i,j \leq k-1} u_{i,j} \geq (k-2 - \varepsilon) \log X$ replaced by the condition that $\sum_{i,j \leq k-1} u_{i,j} \geq (k-2 + \varepsilon) \log X$). 
It follows that on this set of $r_i$ and $c_i$ we have
\begin{align*}
\text{Vol}({\widetilde{\mathcal{U}}_{k,\varepsilon}}(\mathbf{r}, \mathbf{c}, X))
> \text{Vol} & ({\widetilde{\mathcal{U}}_{k, \varepsilon}}(\mathbf{(1 - \varepsilon) \log X}, \mathbf{(1 - \varepsilon) \log X}, X)) \\
& = (1 + o(1)) (\log X)^{(k-1)^2} \cdot k^{-(k-1)} \text{Vol}(\mathcal{B}_k)
\end{align*}
Therefore we obtained the following lower bound
\begin{align*}
(1 + o(1)) X (\log X)^{(k-1)^2}&  \cdot k^{-(k-1)}  \text{Vol}(\mathcal{B}_k)
\\ & \times \iint_{\mathcal{\widetilde{T}}_{k,\varepsilon}(X)} \exp \Big ( \min ( \sum_{1 \leq j \leq k-1} c_j, 
\sum_{1 \leq i \leq k-1} r_i) \Big ) \prod_{i \leq k-1} dc_i \ dr_i
\end{align*}
where $\widetilde{\mathcal{T}}_{k,\varepsilon}(X)$ is the set of $r_i, c_i$
satisfying $(1 - \varepsilon) \log X < r_i, c_i \leq \log X$
for all $i \leq k-1$. Note that the condition $(1 - \varepsilon) \log X<r_i, c_i$
can be dropped. Indeed the contribution to the integral 
of any tuple of $(r_1, \ldots, r_{k-1})$ or $(c_1,\ldots, c_{k-1})$
where at least one of the $c_i, r_i$ is $\leq (1 - \varepsilon) \log X$ is $\leq X^{k-1 - \varepsilon}$ and therefore negligible. Thus we
can extend the integration to all of $c_i, r_i \leq \log X$. Because of this
 Lemma \ref{lem:mainlemma} is applicable and we have therefore obtained the lower bound
$$
\geq (1 + o(1)) \binom{2k-2}{k-1} \cdot k^{-(k-1)} \text{Vol}(\mathcal{B}_k)
X^k \cdot (\log X)^{(k-1)^2}
$$
as claimed. 
Since we have obtained asymptotically matching upper and lower bounds the proof of the proposition is finished. 
\end{proof}

\section{Proof of Theorem 4}
In the Rademacher case we have, letting $\square$ denote a generic square,
\begin{align*}
\mathbb{E} \Big ( \sum_{n \leq X} f(n) \Big )^{k} & = \sum_{\substack{n_1, \ldots, n_k \leq X}}
\mathbb{E}[f(n_1) \ldots f(n_k)] \\
& = \sum_{\substack{n_1, \ldots, n_k \leq X \\
n_1 \ldots n_k = \square}} \mu^2(n_1) \ldots \mu^2(n_k)
\end{align*}
Let $g(n_1, \ldots, n_k)$ be a multiplicative
function of several variables, supported on square-free $n_i$, and such
that $g = 1$ when $n_1 \ldots n_k = \square$ and $g = 0$ otherwise. 
Then we find that the Dirichlet series
$$
\sum_{n_1 = 1}^{\infty} \ldots \sum_{n_k = 1}^{\infty}
\frac{g(n_1, \ldots, n_k)}{n_1^{s_1} \ldots n_k^{s_k}}
$$
is equal to
$$
\prod_{p} \Big ( 1 + \sum_{\substack{0 \leq \alpha_1, \ldots, \alpha_{k} \leq 1
\\ \alpha_1 + \ldots + \alpha_{k} \equiv
0 \mod{2}}} \frac{1}{p^{\alpha_1 s_1 + \ldots + \alpha_k s_k}} \Big ) .
$$
This factors as
$$
H(s_1, \ldots, s_{k}) \prod_{1 \leq i < j \leq k} \zeta(s_i + s_j) 
$$
with
$$
H(\tfrac 12, \ldots, \tfrac 12) = \prod_{p} 
\Big ( 1 - \frac{1}{p} \Big )^{k(k-1)/2}
\Big (1 + \sum_{1 \leq j \leq k/2}
\frac{\binom{k}{2j}}{p^j} \Big )
$$
The main result of La Bret\`eche is applicable with $\binom{k}{2}$ linear forms
$\ell^{(i,j)}(s_1, s_2, \ldots, s_k) = s_i + s_j$ defined for $1 \leq i < j \leq k$. The rank of these
linear forms is equal to $k$ for $k \geq 3$ (for $k = 2$ the rank is equal
to $1$ since there is only one form in that case). Therefore applying
La Bret\`eche's result it follows that the moment is asymptotically
$$
(1 + o(1)) C_k X^{k/2} (\log X)^{\binom{k}{2} - k}. 
$$
In order to determine the constant $C_k$ one could use Th\'eor\`eme 2 of La Bret\`eche, to
conclude that the moment is asymptotically
$$
(1 + o(1)) H(\tfrac 12, \ldots, \tfrac 12) \text{Vol}(B(X))
$$
where $B(X)$ is the set of $(u_{i,j})_{i < j} \in \mathbb{R}^{k(k-1)/2}$ such that
\begin{align*}
\text{ for all } 1 \leq i \leq k: \ \prod_{j < i} u_{j,i} \prod_{i < j} u_{i,j}  
\leq X
\end{align*}
and then proceed in a manner similar to Theorem 3. However we leave this computation to the interested reader.


\end{document}